\documentclass{amsart}
\usepackage{}
\usepackage{color}
\usepackage{amssymb}

\newtheorem{theorem}{Theorem}[section]

 \theoremstyle{definition}

\theoremstyle{remark}

\numberwithin{equation}{section}

\begin{document}

\title[Neumann problem for Hessian equation]
{Interior and boundary gradient estimates for Neumann problem of fully nonlinear Hessian equations }

\author{Weisong Dong}
\address{School of Mathematics, Tianjin University,
         Tianjin, 300354, China}
\email{dweeson@gmail.com}

%\date{}

\begin{abstract}

In this paper we study the {\it a priori} gradient estimates for admissible solutions to Neumann boundary value problem
of fully nonlinear Hessian equations on Riemannian manifolds. We firstly derive an interior gradient estimates for admissible
solutions, then we obtain boundary gradient estimates based on the interior gradient estimates we have got.



\end{abstract}

\maketitle

\section{Introduction}

Let $(M,g)$ be a smooth compact Riemannian manifold of dimension $n\geq 2$ with Riemannian metric $g$ and
nonempty boundary $\partial M$. We denote the interior of $M$ by $\mathring M$, i.e. $\mathring M = M/\partial M$.
In this paper, we firstly study the {\it a priori} interior gradient estimates for solutions to the fully nonlinear equation
\begin{equation}
\label{eqn}
 f (\lambda(\nabla^2 u))  =  \psi (x,u, \nabla u),\; \mbox{on $M$},
\end{equation}
where $\nabla$ is the Levi-Civita connection of $M$, $\nabla u$ is the gradient of function $u$,
$\lambda (\nabla^2 u)$ denotes the eigenvalues of
the Hessian $\nabla^2 u$ of function $u$ with respect to $g$,
and $\psi\in C^3 (M\times\mathbb{R}\times\mathbb{R}^n) >0$.
The smooth function $f$ defined on an open convex symmetric cone
$\Gamma\subset\mathbb{R}^n$ with vertex at the origin is symmetric with respect to $\lambda\in \Gamma$.
The cone also satisfies the following condition
\[ \Gamma_n \equiv\{\lambda\in\mathbb{R}^{n}:
\mbox{each component } \lambda_{i}>0\} \subseteq \Gamma \neq \mathbb{R}^{n}. \]
We call a function $u$ is an admissible solution of equation \eqref{eqn} if
$\lambda (\nabla^2 u)\in \Gamma$ holds on $M$. In this work,
we derive gradient estimates for admissible solutions.

For the function $f$, we suppose it
satisfies the following conditions:
\begin{equation}
\label{f1}
f_i\equiv \frac{\partial f}{\partial \lambda_i} >0 ,
\end{equation}
\begin{equation}
\label{f2}
f \mbox{\;is concave},
\end{equation}

\begin{equation}
\label{f3}
f > 0 \;\mbox{in}\; \Gamma, f = 0 \;\mbox{on} \;\partial \Gamma,
\end{equation}
and
\begin{equation}
\label{f5}
f \,\mbox{is homogeneous of degree one}.
\end{equation}
These conditions are well known in the study of fully nonlinear elliptic equations
since the fundamental work in \cite{CNS}.
We mention some results in \cite{CNS} here.
Let $F(r) = f(\lambda(r))$ for $r\in \mathcal{S}^{n\times n}$ with
$\lambda(r)\in \Gamma$, where $\mathcal{S}^{n\times n}$ is the set
of $n\times n$ symmetric matrices.
Denote $F^{ij}:=\frac{\partial F}{\partial r_{ij}}$.
The first result is that the matrix $\{F^{ij}\}$ has eigenvalues $f_1, \cdots, f_n$.
By \eqref{f1},
$\{F^{ij}\}$ is positive definite for $\lambda(r)\in \Gamma$,
which means equation \eqref{eqn} is elliptic for admissible solutions.
Another result is that the assumption \eqref{f2} implies that
$F$ is concave.
To derive both interior and boundary gradient estimates, we have to assume
the following condition
\begin{equation}
\label{f6}
f_j (\lambda) \geq \nu_0 \Big(1 + \sum f_i (\lambda)\Big), \mbox{ for any }
  \lambda \in \Gamma \mbox{ with } \lambda_j < 0,
\end{equation}
where $\nu_0$ is a uniform positive constant. It is wildly used in deriving gradient estimates, see
\cite{GS}, \cite{T} and \cite{U}.
We also assume that the following growth condition for $\psi(x,z,p)$
holds in $M\times\mathbb{R}\times\mathbb{R}^n$:
\begin{equation}
\label{psi}
|\psi_{x_i}| + |\psi_{z}||\nabla u| + |\psi_{p_l}| |\nabla u|^2\leq C |\nabla u|^{2+\gamma}
\end{equation}
 for a positive constant $\gamma<1$ and a uniform positive constant $C$.

The classical examples of $f$ are given by
$(\sigma_k / \sigma_l)^{1/(k - l)}$,
$0 \leq l < k \leq n$, defined on the cone
$\Gamma_{k} = \{\lambda \in \mathbb{R}^{n}: \sigma_{j} (\lambda) > 0, j = 1, \ldots, k\}$,
where $\sigma_{k} (\lambda)$ is the $k$-th elementary symmetric function
\[
\sigma_{k} (\lambda) = \sum_ {i_{1} < \ldots < i_{k}}
\lambda_{i_{1}} \ldots \lambda_{i_{k}},\ \ k = 1, \ldots, n,
\]
and we define that $\sigma_0 \equiv 1$. In \eqref{eqn}, if $f=\sigma_n^{1/n}$, it is
the well known Monge-Amp\`ere equation. When $f = (\sigma_k)^{1/k}$, it is called $k$-Hessian
equation.

Throughout this paper, we suppose the assumptions \eqref{f1}-\eqref{psi} hold.
Our first result is the following interior
gradient estimates for admissible solutions of equation \eqref{eqn} on Riemannian manifold $(M, g)$.

\begin{theorem}
\label{main1}
Suppose $u\in C^3(M)$ is an admissible solution of equation \eqref{eqn}. For any geodesic ball $B_r (x_0)\subset \mathring M$ with center $x_0\in M$ and radius
$r<1$, there exists a positive
constant $C$ such that
\[
|\nabla u(x_0)| \leq \frac{C}{r},
\]
where $C$ only depends on $|u|_{C^0}$ and other known data.
\end{theorem}

The next result of this work is the {\it a priori} boundary gradient estimates for admissible solutions
to equation \eqref{eqn} with the Neumann boundary value condition
\begin{equation}
\label{eqn-b}
u_{\nu} = \varphi (x, u), \; \mbox{on $\partial M$},
\end{equation}
where $\nu$ is the unit inner normal vector to $\partial M$ with respect to $g$
and the function $\varphi \in C^{3}(\partial M\times \mathbb{R})$.
In the recent breakthrough work \cite{MX} and \cite{MQ} for Neumann problem,
the boundary gradient estimates are established by reduction
to the interior gradient estimates.
Here, we follow their idea. Given a small positive constant $\mu$, we define the boundary strip of $\partial M$ by
\[
M_\mu := \{x\in M, 0< d (x, \partial M) <\mu\}.
\]
We prove that the gradient of admissible solutions can be bounded uniformly on $\overline{M_\mu}$
when $\mu$ is sufficiently small, where the bound depends on the interior gradient estimates.
%This method has been used in many literatures concerning Neumann type problems, such as \cite{S},  \cite{JXX} and \cite{GX}.

Now we state the boundary gradient estimates for admissible solutions of
equation \eqref{eqn} with Neumann boundary value \eqref{eqn-b}.
\begin{theorem}
\label{main2}
Suppose $u\in C^3(M)$ is an admissible solution of \eqref{eqn} with Neumann boundary value \eqref{eqn-b}.
Then, there is a small positive constant $\mu$ such that
\[
\sup_{M_\mu} |\nabla u| \leq C (1 + \sup_{M/M_\mu} |\nabla u|) ,
\]
where $C$ depends on $|u|_{C^0}$, $\mu$ and other known data.
\end{theorem}

Combining Theorem 1.1 and Theorem 1.2 we can get the following gradient estimates for admissible solutions of
Neumann problem \eqref{eqn} and \eqref{eqn-b}.
\begin{theorem}
\label{main}
Suppose $u\in C^3(M)$ is an admissible solution of Neumann problem \eqref{eqn} and \eqref{eqn-b}. We have
\[
\sup_M |\nabla u| \leq C,
\]
where $C$ depends on $|u|_{C^0}$ and other known data.
\end{theorem}

Theorem 1.1 extends some previous work by Chou and Wang \cite{CW} and Chen \cite{Chen}, in which
the interior gradient estimates for admissible solutions are carried out for equation \eqref{eqn}
with $f=\sigma^{1/k}_k$ and $f=(\sigma_k / \sigma_l)^{1/(k - l)}$
respectively in Euclid space.
%On closed compact
%Riemaniann manifold with non-negative sectional curvature,
%Li \cite{Li} get the global gradient estimates for the special type
%$f(\lambda(\nabla^2 u+ \mbox{Id}))=\psi$.
Recently, the interior gradient estimates for more general Monge-Amp\`ere type equation have been derived
in Euclid space \cite{JXX} and
on Riemaniann manifolds \cite{GX} without the growth condition on $\psi$.

In conformal geometry, there are interior gradient estimates for fully nonlinear equations,
see \cite{GW}, \cite{GLW} and \cite{JLL}.
For the prescribed mean curvature equation, such estimates have been extensively studied, see \cite{GT}.
Interior gradient estimates have also been obtained for Weingarten curvature equations in \cite{Ko}
and see \cite{T} for some other high order curvature equations. Li \cite{Li-JDE} generalized the results in \cite{Ko} to
general curvature equations
by analysis on the graphs of solutions.
We remark that Wang in \cite{W} obtained the interior gradient estimates for Weingarten curvature
equations using the maximum principle with suitable choice of auxiliary function.
For the interior gradient estimates of general curvature equations as in \cite{Li-JDE}, there is no such proof.

It is well known that the {\it a priori} estimates are crucial to
the proof of existence of solutions for Neumann problem \eqref{eqn} and \eqref{eqn-b} by the continuity method.
We note that the existence result of the Neumann problem for equations of Monge-Amp\`ere type
was obtained in \cite{LTU}.
Recently, the Neumann problem for $k$-Hessian equation with $\varphi_z <0$ is solved by Ma and Qiu \cite{MQ},
and see Qiu and Xia \cite{QX} for the case $\varphi_z \equiv 0$.
Ma and Xu \cite{MX} have got the boundary gradient estimates for the
Neumann boundary value problem of mean curvature equation.
The prescribed contact angle problem for the mean curvature equation was studied a long time ago, see
\cite{Ural},\cite{S},\cite{SS},\cite{Gerh}.

The rest of the paper is organized as follows: in Section 2, we prove the interior gradient estimates
for fully nonlinear equation \eqref{eqn}.
In section 3, we derive gradient estimates on $\overline{M_\mu}$ for equation \eqref{eqn} with Neumann boundary
value \eqref{eqn-b}  via the interior gradient estimates.

\section{Proof of Theorem 1.1}

Without loss of generality, we suppose $0\in\mathring{M}$.
Choose a positive constant $r<1$ such that $B_r(0)\subset \mathring{M}$.
We prove there exists a uniform constant $C$ depending on $|u|_{C^0}$ such that $|\nabla u(0)|\leq \frac{C}{r}$.
Consider the following auxiliary function on $B_r(0)$,
\[
G(x) = |\nabla u| h (u) \zeta(x),
\]
where $\zeta(x) = r^2-\rho^2(x)$ and $\rho(x)$ denotes the geodesic distance from $0$.
We can assume $r$ is small such that $\rho$ is smooth and $|\nabla\rho|\equiv 1$ in $B_r(0)$.
Note that $\nabla_{ij} \rho^2(0) = 2\delta_{ij}$. For sufficiently small $r$ , we may further assume in $B_r(0)$
that
\[
\delta_{ij}\leq \nabla_{ij} \rho^2 \leq 3\delta_{ij}.
\]
Then we assume that $G(x)$ attains its
maximum at some point $x_0\in B_r(0)$. Choose smooth orthonormal local
frames~$e_1,\cdots, e_n$ about $x_0$ such that
\[
\nabla_1 u(x_0) = |\nabla u(x_0)| ,
\]
and $\{\nabla_{ij} u\}_{2\leq i,j\leq n}$ is diagonal at $x_0$.
We may also assume $\nabla_{i}e_j = 0$ at $x_0$.
We differentiate $\log G$ at $x_0$ twice to get that
\begin{equation}
\label{dif1}
\frac{\nabla_{i1} u}{\nabla_1 u} + \frac{\nabla_i h}{h} + \frac{\nabla_i\zeta}{\zeta}=0
\end{equation}
and
\begin{equation}\label{dif2}\begin{aligned}
0\geq&\frac{\nabla_{ij1}u}{\nabla_1 u} + \frac{\nabla_{ik}u\nabla_{jk}u}{(\nabla_1 u)^2}- \frac{2\nabla_{i1} u\nabla_{j1}u}{(\nabla_1 u)^2} \\
&+ \frac{\nabla_{ij}h}{h}-
\frac{\nabla_i h \nabla_j h}{h^2} +\frac{\nabla_{ij}\zeta}{\zeta}-\frac{\nabla_i\zeta\nabla_j\zeta}{\zeta^2}.
\end{aligned}\end{equation}
Recall the formula for interchanging order of covariant derivatives
\[\nabla_{ijk}u=\nabla_{kij}u + R^l_{kij}\nabla_l u.\]
We have
\begin{equation}\label{deqn}\begin{aligned}
F^{ij}\nabla_{ij1}u
%=&\, F^{ij}(\nabla_{1ij}u+ R^l_{ji1}\nabla_l u)\\
=&\, \psi_{x_1} + \psi_z\nabla_1 u + \psi_{p_l}\nabla_{1l}u + F^{ij} R^l_{ji1}\nabla_l u.\\
\end{aligned}\end{equation}
%where $\nabla'_k$ denote the covariant derivatives with respect to $x_k$.
By the growth condition of $\psi$, we have
\begin{equation}
\label{deqn'}
F^{ij}\nabla_{ij1}u \geq  -C |\nabla_1 u|^{2+\gamma} \Big(1+ \sum F^{ii}\Big)+ \psi_{p_l}\nabla_{1l}u,
\end{equation}
where $C$ is a uniform positive constant.
By equality \eqref{dif1} we see that
\[\begin{aligned}
F^{ij}\frac{\nabla_{i1} u\nabla_{j1}u}{(\nabla_1 u)^2}
=&\,F^{ij}\Big(\frac{\nabla_i h}{h} + \frac{\nabla_i\zeta}{\zeta}\Big)
\Big(\frac{\nabla_j h}{h} + \frac{\nabla_j\zeta}{\zeta}\Big)\\
\leq&\,2F^{ij}\frac{\nabla_i h}{h}\frac{\nabla_j h}{h}
+ 2F^{ij} \frac{\nabla_i\zeta}{\zeta}\frac{\nabla_j\zeta}{\zeta}.
\end{aligned}\]
It is readily to see
\[
F^{ij}\frac{\nabla_{ij}\zeta}{\zeta}-F^{ij}\frac{\nabla_i\zeta\nabla_j\zeta}{\zeta^2}\geq -\frac{3\sum F^{ii}}{\zeta} - \frac{4\rho^2F^{ij}\nabla_i\rho\nabla_j \rho}{\zeta^2}.
\]
Since $\{F^{ij}\}> 0$,
contracting \eqref{dif2} with $\{F^{ij}\}$ and by \eqref{deqn'} and \eqref{dif1} we get
\begin{equation}\label{sum}\begin{aligned}
0\geq &\, -C |\nabla_1 u|^{1+\gamma}\Big(1+\sum F^{ii}\Big) - \frac{h'}{h} \psi_{p_1}\nabla_1 u + \frac{2\rho}{\zeta}\psi_{p_l}\nabla_l \rho+\frac{h'}{h}F^{ij}\nabla_{ij}u \\
&\,+ \Big(\frac{h''}{h}-3\frac{h'^2}{h^2}\Big)F^{11}(\nabla_1 u)^2 -\frac{3}{\zeta}\sum F^{ii} - \frac{12r^2}{\zeta^2}F^{ij}\nabla_i\rho\nabla_j \rho,
\end{aligned}\end{equation}
here and in what follows $C$ will be a positive constant only depending on $|u|_{C^0}$
and other known data, but may change from line to line.

Now we determine the function $h$. Let $\delta$ be a small positive constant to be chosen and
$C_0$ be a positive constant such that $|u|_{C^0}+ 1\leq C_0$. Let $h$ be defined by
\[
h(u) = e^{\delta (u+C_0)^2}.
\]
Differentiating $h$ twice we get
\[
h' = 2\delta (u + C_0 ) e^{\delta(u+C_0)^2} = 2\delta (u + C_0 ) h
\]
and
\[
h'' = 4\delta^2 (u + C_0)^2 h + 2\delta h.
\]
It follows that
\[
\frac{h''}{h}-3\frac{h'^2}{h^2} = 2\delta - 8\delta^2 (u+C_0)^2 \geq \delta
\]
when $\delta$ is sufficiently small.

By the homogeneity of $f$, we see that $F^{ij}\nabla_{ij}u = \psi \geq 0 $. With the assumption $f$ is homogeneity and
$f$ is concave, we have
\[
\sum f_i (\lambda) = f(\lambda) + \sum f_i (\lambda) (1-\lambda_i) \geq f(\mathbf{1})
\]
where $\mathbf{1}= (1, \cdots, 1)\in \Gamma$. Without loss of generality, we may assume that $f$ is normalized such that
$ f(\mathbf{1})=1$. We now derive from \eqref{sum} that
\begin{equation}\label{sum2}\begin{aligned}
0\geq &\,  \delta F^{11}(\nabla_1 u)^2- C|\nabla_1 u|^{1+\gamma}\Big(1+\sum F^{ii}\Big)
 -\Big(\frac{3}{\zeta} + \frac{Cr\nabla_1 u}{\zeta}+\frac{12r^2}{\zeta^2}\Big)\sum F^{ii}
\end{aligned}\end{equation}
where we used the fact $|\nabla \rho|\equiv 1$ and $|\psi_{p_l}| \leq C|\nabla u|$ when $|\nabla u|$ is
sufficiently large.
From \eqref{dif1} we have
\[\begin{aligned}
\nabla_{11}u =&\, -\nabla_1 u \Big(\frac{h'}{h}\nabla_1 u - \frac{2\rho}{\zeta}\nabla_1\rho\Big)\\
=&\,-2\delta (u+C_0)(\nabla_1u)^2 +  \frac{2\rho}{\zeta}\nabla_1 u\nabla_1\rho.
\end{aligned}\]
Since $u+C_0\geq 1$,
we assume at $x_0$ that $\nabla_1 u$ is sufficiently large such that $\zeta\nabla_1 u > \frac{\rho\nabla_1 \rho}{\delta (u+C_0)}$.
This means that $\nabla_{11} u< 0$ at $x_0$.
By the assumption \eqref{f6} we now have
\[
F^{11} \geq \nu_0 + \nu_0 \sum F^{ii}.
\]
Hence, when $\nabla_1 u(x_0)$ is sufficiently large such that
$|\nabla_1 u(x_0)|^{1-\gamma}\geq \frac{C}{\delta \nu_0}$, from \eqref{sum2} we have
\[\begin{aligned}
0\geq &\, \frac{\delta}{2} F^{11}(\nabla_1 u)^2
 -\Big(\frac{3}{\zeta} + \frac{Cr\nabla_1 u}{\zeta}+\frac{12r^2}{\zeta^2}\Big)\sum F^{ii}\\
 \geq &\, \frac{\sum F^{ii}}{\zeta^2} \Big(\frac{\delta\nu_0}{2}(\zeta\nabla_1 u)^2  -Cr\nabla_1 u\zeta-15r^2\Big).
\end{aligned}\]
We therefore obtain that
\[
\zeta(x_0)\nabla_1 u(x_0) \leq  \frac{Cr}{\delta\nu_0}.
\]
We just proved that $G(x_0)$ is bounded by some uniform positive constant $C$, which depends only on $|u|_{C^0}$
and other known data. Therefore, by
$G(0) \leq G(x_0)$, we have
\[
|\nabla u(0)| \leq \frac{C }{r^2} \zeta(x_0)\nabla_1 u(x_0) \leq \frac{C}{r},
\]
where $C$ depends only on $|u|_{C^0}$
and other known data.
The proof of Theorem 1.1 is complete.

\section{Proof of Theorem 1.2}

%We assume that $\mu$ is sufficiently small such that $d(x)|\varphi|_{C^3} \leq \frac{1}{2}$ in $M_\mu$.
In this section, we suppose that $\varphi$ is smoothly extended to $M\times \mathbb{R}$.
In order to deal with the boundary gradient estimates for Neumann problem,
the following function
\[
w(x):= u(x) - \varphi (x, u) d(x)
\]
has been used in the construction of auxiliary function, which
ensures that the maximum point of the auxiliary function must be an interior point of $M$.
This function has been used in recent pioneering work on Neumann problem
of $k$-Hessian equation, see \cite{MQ} and \cite{QX}, and of the mean curvature equation, see \cite{MX}.

Now we consider the auxiliary function
\[
G(x):= |\nabla w|^2 e^{h (u)} e^{\phi (d)}
\]
where
$h(u)= \delta(u+C_0)^2$ and
$\phi(d)= A d$, here $C_0$ is as in Section 2 and $A$ is a positive constant to be chosen.
Assume $\max_{x\in\bar M_{\mu}} G(x)$ is attained at $x_0$. We prove Theorem 1.2 by three cases.

For the first case that $x_0 \in \partial M$, we choose local orthonormal coordinate at $x_0$ such that
$e_n = \nu$ and $e_1, \cdots, e_{n-1}$ are tangential to the boundary $\partial M$. We have at $x_0$ that
\begin{equation}
\label{bd1}
0 \geq ({\log G})_\nu = \frac{\nabla_\nu |\nabla w|^2}{|\nabla w|^2} + h' \nabla_\nu u + A \nabla_\nu d.
\end{equation}
Note that $\nabla_\nu d = 1$ and $d = 0$ on $\partial M$.
%Since $\nabla w= \nabla u - \nabla \varphi d - \varphi \nabla d$ and $|\nabla d|\equiv 1$, we have
%\[
%|\nabla w|^2 = |\nabla u|^2 + |\nabla\varphi|^2d^2 + \varphi^2 - 2 d g(\nabla u, \nabla \varphi)
%- 2\varphi g(\nabla u, \nabla d) + 2d\varphi g(\nabla \varphi, \nabla d)
%\]
%By the boundary condition $g(\nabla u, \nabla d) =\varphi$, we have
%\[\begin{aligned}
%\nabla_\nu |\nabla w|^2 =  \nabla_\nu |\nabla u|^2 + \nabla_\nu \varphi^2   - 2 g(\nabla u, \nabla \varphi)
%- 2 \varphi\nabla_\nu g(\nabla u, \nabla d)\\
%& + 2\varphi g(\nabla \varphi, \nabla d)- 2\nabla_\nu \varphi g(\nabla u, \nabla d)
%\end{aligned}\]
%Under an orthonormal frame $e_1,\cdots, e_{n-1},e_n=\nu$ around $x_0$, we have
%\[\begin{aligned}
%&\nabla_\nu |\nabla u|^2 = 2 \sum_{i=1}^n\nabla_i u\nabla_{ni}u\\
%&\nabla_\nu \varphi^2 = 2\varphi  \nabla_\nu\varphi = 2 \nabla_\nu u \nabla_{\nu} \varphi\\
%&\nabla_\nu g(\nabla u, \nabla d) = \nabla_{n}\varphi
%\end{aligned}\]
%\[
%\nabla_\nu |\nabla w|^2 = 2 \sum_{i=1}^n\nabla_i u\nabla_{ni}u - 2\sum_{i=1}^{n-1} \nabla_i u\nabla_i \varphi -2\varphi\nabla_{n}\nabla_n\varphi
%\]
We see that
\[
\nabla_\nu w = \nabla_\nu u - \varphi \nabla_\nu d = 0.
\]
Then \eqref{bd1} becomes to
\[
0 \geq \frac{2 \nabla_\alpha w \nabla_{\nu \alpha} w}{|\nabla w|^2} + h' \nabla_\nu u + A,
\]
where
%we use the Einstein summation convention and
repeated $\alpha$ means summation from $1$ to $n-1$.
With the boundary condition \eqref{eqn-b} we have
\[\begin{aligned}
\nabla_{\nu \alpha } w = &\, \nabla_{\alpha\nu} u - \nabla_{\alpha\nu} (\varphi d)\\
= &\, \nabla_{\alpha\nu} u -  d  \nabla_{\alpha\nu} \varphi- \varphi \nabla_{\alpha \nu}d - \nabla_\alpha \varphi \nabla_{\nu} d
 -\nabla_\alpha d \nabla_{\nu} \varphi \\
  = &\, \nabla_\alpha \varphi -\nabla_{\nabla_\alpha \nu} u -\varphi \nabla_{\alpha \nu}d
  - \nabla_{\alpha} \varphi \nabla_{\nu} d- \nabla_\alpha d \nabla_{\nu} \varphi.
\end{aligned}\]
We should note that $|\nabla w|$ and $|\nabla u|$ are equivalent when $|\nabla u|$ is large enough. It yields that on $\partial M$
\[
|\nabla_{\nu \alpha } w | \leq C |\nabla w|,
\]
where $C$ only depends on $|u|_{C^0}$ and other known data.
Hence, we have
\[
0 \geq - C - h'|\varphi|_{C^0} + A,
\]
from which we can get a contradiction if $A$ is sufficiently large.
We therefore have that $G(x_0)$ is bounded.

For the second case that $x_0 \in M_\mu$, we use the maximum principle. Differentiating $\log G$ at $x_0$ twice we obtain
\begin{equation}
\label{II-1}
\frac{2 \nabla_k w \nabla_{jk} w}{|\nabla w|^2} + h' \nabla_j u + A \nabla_j d = 0,
\end{equation}
and
\begin{equation}
\label{II-2}
\begin{aligned}
0 \geq &\,\frac{1}{|\nabla w|^2}\{2 \nabla_k w \nabla_{ijk}w + 2 \nabla_{ik}w \nabla_{jk} w \\
&\,- \frac{2 }{|\nabla w|^2}  \nabla_k w \nabla_l w \nabla_{jk}w \nabla_{il} w\}\\
 &\,+ h'' \nabla_i u \nabla_j u + h' \nabla_{ij} u + A \nabla_{ij} d.
\end{aligned}
\end{equation}
Contracting \eqref{II-2} with $\{F^{ij}\}$, by \eqref{II-1} and $\{F^{ij}\}>0$, we have
\begin{equation}
\label{II-4}
\begin{aligned}
0 \geq &\, \frac{2}{|\nabla w|^2} F^{ij } \nabla_k w \nabla_{ijk}w +(h''- 2h'^2)  F^{ij}\nabla_i u\nabla_j u \\
&\, - 2A^2 F^{ij}\nabla_i d\nabla_j d + h' F^{ij}\nabla_{ij} u + A F^{ij}\nabla_{ij} d.
\end{aligned}\end{equation}
For convenience, we denote that $\Phi(x,u) = \varphi(x,u) d(x)$. Then, $w = u - \Phi $.
By direct calculating we see
\begin{equation}
\label{dw}
\nabla_i w =(1 - \Phi_u) \nabla_i u - \Phi_{x_i}.
\end{equation}
We can choose $\mu$ sufficiently small such that $\frac{1}{2}\leq1-\Phi_u\leq 1$.
Then, $|\nabla w|$ and $|\nabla u|$ is equivalent in $M_\mu$
when $|\nabla u|$ is sufficiently large.
Differentiating $\nabla_k w$ again, we have
\[\begin{aligned}
\nabla_{ik}w = (1-\Phi_u)\nabla_{ik}u - (\Phi_{ux_i} + \Phi_{uu}\nabla_i u)\nabla_k u - \Phi_{x_k u} \nabla_i u - \Phi_{x_kx_i}.
\end{aligned}\]
Since $\Phi_{uu} = \varphi_{uu} d$, when $|\nabla u|$ is sufficiently large, we have
\[
\nabla_{ik} w  = (1 - \Phi_u) \nabla_{ik} u + d O( |\nabla w|^2) + O (|\nabla w|)
\]
and
\begin{equation}
\label{II-3}
\begin{aligned}
\nabla_k w \nabla_{ik}w = &\,  (1- \Phi_u) \nabla_k w  \nabla_{ik} u + d O( |\nabla w|^3) + O (|\nabla w|^2).
\end{aligned}
\end{equation}
From \eqref{II-1} we see, for any $1\leq i\leq n$,
\begin{equation}
\label{key}
\nabla_k w \nabla_{ik} u =  O(|\nabla u|^3).
\end{equation}
Differentiating $\nabla_{jk}w$ again, we see that
\[\begin{aligned}
\nabla_{ijk} w = &\, (1- \Phi_u) \nabla_{ijk} u - (\Phi_{ux_i}+ \Phi_{uu}\nabla_i u)\nabla_{jk}u
 - (\Phi_{ux_j} + \Phi_{uu}\nabla_j u)\nabla_{ik} u \\
&\, - (\Phi_{ux_jx_i} + \Phi_{ux_j u} \nabla_i u+ \Phi_{uux_i}\nabla_j u + \Phi_{uuu}\nabla_iu\nabla_ju + \Phi_{uu}\nabla_{ij}u)\nabla_k u\\
&\, - \Phi_{x_k u}\nabla_{ij}u - (\Phi_{x_k u x_i}+ \Phi_{x_k u u} \nabla_i u )\nabla_j u - \Phi_{x_k x_j x_i} - \Phi_{x_kx_j u}\nabla_i u.
\end{aligned}\]
Differentiating the equation \eqref{eqn} and by the growth condition of $\psi$, we get
\[
F^{ij}\nabla_{kij}u = \psi_{x_k} + \psi_z \nabla_k u + \psi_{p_l} \nabla_{kl} u = O(|\nabla u|^{2+\gamma}) + O(|\nabla u|^{\gamma}) \nabla_{kl}u,
\]
when $|\nabla u|$ is sufficiently large.
%\[
%\nabla_k A^{ij} = \nabla_k' A^{ij} + A^{ij}_z \nabla_k u + A^{ij}_{p_l} \nabla_{kl} u = O(|\nabla u|^{2+ \gamma})  + O(|\nabla u|^{\gamma}) \nabla_{kl} u
%\]
Recall the formula for interchanging order of covariant derivatives
\[\nabla_{ijk}u=\nabla_{kij}u + R^l_{kij}\nabla_l u.\]
We have
\[\begin{aligned}
F^{ij} \nabla_{ijk} u %= &\, F^{ij} \nabla_{kij} u+ F^{ij}R^l_{ijk} \nabla_l u \\
%=&\, F^{ij} \big(\nabla_k U_{ij} -\nabla_k A^{ij}\big) + O(|\nabla u|)\sum F^{ii}\\
= O(|\nabla u|^{2+ \gamma}) + O(|\nabla u|^{\gamma}) \nabla_{kl}u  + O(|\nabla u|)\sum F^{ii}.
\end{aligned}\]
Note that $\Phi_{uuu} = \varphi_{uuu} d$. We have
\[\begin{aligned}
F^{ij} \nabla_{ijk} w = &\, (1- \Phi_u) F^{ij} \nabla_{ijk} u - 2 F^{ij} \big(\Phi_{ux_i}+ \Phi_{uu}\nabla_i u \big)\nabla_{jk}u\\
 &\, - F^{ij} \nabla_{ij} u \big(\Phi_{uu}  \nabla_k u + \Phi_{x_ku}\big)+ \big( d O(|\nabla u|^3) + O(|\nabla u|^2) \big) \sum F^{ii}\\
 =&\, O(|\nabla u|^{\gamma}) \nabla_{kl}u + O(|\nabla u|^{2+ \gamma})- 2 F^{ij} \big(\Phi_{ux_i}+ \Phi_{uu}\nabla_i u \big)\nabla_{jk}u\\
 &\, - F^{ij} \nabla_{ij} u \big( \Phi_{uu}  \nabla_k u + \Phi_{x_ku}\big)+ \big( d O(|\nabla u|^3)  + O(|\nabla u|^2) \big) \sum F^{ii}.
\end{aligned}\]
Using \eqref{key} we derive
\[\begin{aligned}
F^{ij} \nabla_k w \nabla_{ijk}w
=&\, O(|\nabla u|^{3+\gamma}) - F^{ij} \nabla_{ij} u O(|\nabla u|) \Big(d |\nabla u|  \\
&\, + 1 \Big)
 +  \Big( d O(|\nabla u|^4) + O(|\nabla u|^{3}) \Big) \sum F^{ii}.
\end{aligned}\]
By \eqref{II-2}, we obtain, when $|\nabla u|$ is sufficiently large, that
\begin{equation}
\label{II-5}
\begin{aligned}
0 \geq &\, \big(h''- 2(h')^2\big)F^{ij}\nabla_i u\nabla_j u + O(|\nabla u|^{1+\gamma})\\
&\,+ \Big(h' - d O(1) - \frac{O(1)}{|\nabla u|}\Big)F^{ij} \nabla_{ij} u \\
&\, + \big( d O(|\nabla u|^2) + O(|\nabla u|) \big) \sum F^{ii}.
\end{aligned}
\end{equation}
By our definition $h(u) = \delta (u + C_0)^2$ and the choice of $C_0$, one can see
\[
h'(u) = 2\delta (u + C_0)\geq 2\delta, \; h''(u) = 2\delta.
\]
We choose $\delta$ small such that
\[
h''-2(h')^2 = 2\delta - 8\delta^2(u + C_0)^2 >\delta.
\]
Without loss of generality we assume that $\{\nabla_{ij} u (x_0)\}$ is diagonal.
Now suppose that $|\nabla u| \leq n |\nabla_1 u|$.
We obtain from \eqref{dw} that
\[
|\nabla_1 w| = O(|\nabla_1 u|) = O(|\nabla u|) = O(|\nabla w|).
\]
From \eqref{II-1} and \eqref{II-3}, we have
\[
(1-\Phi_u) \nabla_k w \nabla_{1k} u = -  h' \nabla_1 u  |\nabla w|^2 - d O(|\nabla w|^3)+ O(|\nabla w|^2).
\]
Since $\{\nabla_{ij} u (x_0)\}$ is diagonal, we obtain
 \[\begin{aligned}
 \nabla_{11} u = - h'|\nabla w|^2\frac{\nabla_1 u}{(1-\Phi_u) \nabla_1 w}
 - \frac{d O(|\nabla w|^3)}{(1-\Phi_u) \nabla_1 w} +  O(|\nabla w|)
 \end{aligned}\]
 which implies that $\nabla_{11}u<0$ when $|\nabla u|$ is sufficiently large and $d<\mu$ is sufficiently small.
From \eqref{f6}, we get
\[F^{ij}\nabla_iu\nabla_j u = F^{ii}(\nabla_i u)^2\geq \nu_0|\nabla_1 u|^2\Big(1+ \sum f_i\Big). \]
Noting that $F^{ij}\nabla_{ij}u\geq 0$,
we can see
\[
\Big(h' - d O(1) - \frac{O(1)}{|\nabla u|}\Big) F^{ij}\nabla_{ij}u \geq 0,
\]
when $|\nabla u|$ is sufficiently large and $\mu$ is sufficiently small.
We now derive from \eqref{II-5} that
\[\begin{aligned}
0 \geq &\, \Big(\nu_0 \delta + d O(1) + \frac{O(1)}{|\nabla u|}  \Big)\sum F^{ii}+\nu_0 \delta + \frac{O(|\nabla u|^{\gamma})}{|\nabla u|},
\end{aligned}\]
when $|\nabla u|$ is sufficiently large and $\mu$ is sufficiently small.
We get a contradiction from the above inequality if $|\nabla u|$ is sufficiently large and $d<\mu$ is sufficiently small. Therefore, $G(x_0)$
is bounded in this case.

For the third case $x_0 \in \partial M_{\mu} \backslash \partial M$, a bound for $|\nabla u (x_0)|$ can be obtained by Theorem 1.1.
Therefore, $G(x_0)$ is bounded directly.

Now, we have proved that $\max_{\overline M_\mu}G \leq C (1+\sup_{M/M_\mu}|\nabla u|)$,
where $C$ depends on $|u|_{C^0}$ and other known data.
Consequently, we have proved that for admissible solutions $u$ of equation \eqref{eqn}
with Neumann boundary value condition \eqref{eqn-b}, there exists a uniform positive constant $C$ such that
\[
\sup_{M_\mu} |\nabla u| \leq C (1+\sup_{M/M_\mu}|\nabla u|)
\]
holds, where $C$ depends on $|u|_{C^0}$ and other known data.
Theorem 1.2 is proved.

\end{document}